\documentclass[twoside,11pt]{article}
\title{} \author{} \date{}
\usepackage{latexsym,amssymb,times}
\input amssym.def
\newtheorem{te}{Theorem}[section]

\newtheorem{fac}[te]{Fact}
\newtheorem{cla}[te]{Claim}

\newtheorem{df}[te]{Definition}
\newtheorem{rem}[te]{Remark}
\newtheorem{ex}[te]{Example}

\def\dok{\noindent{\bf Proof. }}
\def\kdok{\hfill $\Box$ \par \vspace*{2mm} }
\def\a{\alpha}

\def\f{\varphi}
\def\p{\psi}
\def\o{\omega}
\def\k{\kappa}

\def\r{\rho}
\def\s{\sigma}

\def\t{\tau}

\def\Q{{\mathbb Q}}

\def\N{{\mathbb N}}
\def\X{{\mathbb X}}
\def\Y{{\mathbb Y}}
\def\Z{{\mathbb Z}}

\def\BL{{\mathbb L}}
\def\BH{{\mathbb H}}
\def\BK{{\mathbb K}}
\def\BI{{\mathbb I}}
\def\BF{{\mathbb F}}
\def\BM{{\mathbb M}}

\def\CT{{\mathcal T}}

\def\L{{\mathcal L}}

\def\CF{{\mathcal F}}
\def\c{{\mathfrak{c}}}

\def\la{\langle}
\def\ra{\rangle}

\def\id{\mathop{\mathrm{id}}\nolimits}
\def\otp{\mathop{\rm otp}\nolimits}

\def\Iso{\mathop{\rm Iso}\nolimits}
\def\Aut{\mathop{\rm Aut}\nolimits}

\def\Sym{\mathop{\rm Sym}\nolimits}

\def\Sent{\mathop{\rm Sent}\nolimits}
\def\Form{\mathop{\rm Form}\nolimits}

\def\Age{\mathop{\mathrm{Age}}\nolimits}

\def\ar{\mathop{\rm ar}\nolimits}
\def\Mod{\mathop{\rm Mod}\nolimits}

\def\tp{\mathop{\rm tp}\nolimits}

\def\Th{\mathop{\rm Th}\nolimits}

\def\Pa{\mathop{\rm Pa}\nolimits}
\def\Lit{\mathop{\rm Lit}\nolimits}
\begin{document}
\thispagestyle{plain}
\begin{center}
           {\large \bf \uppercase{Vaught's Conjecture for Monomorphic Theories}}
\end{center}
\begin{center}
{\bf Milo\v s S.\ Kurili\'c}\footnote{Department of Mathematics and Informatics, University of Novi Sad,
                                      Trg Dositeja Obradovi\'ca 4, 21000 Novi Sad, Serbia,
                                      e-mail: milos@dmi.uns.ac.rs}
\end{center}
\begin{abstract}
\noindent
A complete first order theory of a relational signature is called monomorphic iff
all its models are monomorphic (i.e.\ have all the $n$-element substructures isomorphic, for each positive integer $n$).
We show that a complete theory $\CT $ having infinite models is monomorphic iff it has a countable monomorphic model
and confirm the Vaught conjecture for monomorphic theories.
More precisely, we prove that if $\CT$ is a complete monomorphic theory having infinite models, then the number of its non-isomorphic countable models,
$I(\CT ,\o)$, is either equal to $1$ or to $\c$. In addition,
the equality $I(\CT ,\o)=  1$ holds iff some countable model of $\CT$ is simply definable by an $\o$-categorical linear order on its domain.
\end{abstract}
\section{Introduction}\label{S1}
Let $\CT$ be a countable complete theory with infinite models and
$I(\CT, \o )$ the number of non-isomorphic countable models of $\CT$. It is known
that the cardinal $I(\CT, \o )$ can be any positive integer except 2 (Vaught \cite{Vau}); further,
if $\CT$ is the complete theory of a structure with one equivalence relation having arbitrarily large finite equivalence classes, then
$I(\CT, \o )=\o$ and, finally,  $I(\CT, \o )=\c$, if $\CT$ is the complete theory of the linear order $\la \o ,<\ra$.
It is evident that $I(\CT, \o )\leq \c$ and in 1961 Robert Vaught conjectured that  $\o< I(\CT, \o )< \c$ is impossible
(which is trivially true under the CH). In 1970 Morley proved that $I(\CT, \o )>\o _1$ implies $I(\CT, \o )=\c$
(see \cite{Mor}) and, thus, reduced the problem to the cardinal $\o _1$ (as a possible value of $I(\CT, \o )$, which would, together with $\neg$ CH,
produce a counterexample).

The Vaught conjecture was confirmed for several classes of theories; for example,
for theories of linear orders with unary predicates, by Rubin \cite{Rub} in 1974;
for theories of linearly ordered structures with Skolem functions, by Shelah \cite{She} in 1978;
for $\o$-stable theories, by Shelah, Harrington and Makkai \cite{SHM} in 1984;
for o-minimal theories, by Mayer \cite{May} in 1988.
These results were extended by many authors.
In this article we use the following consequence of Rubin's results from \cite{Rub} (see also \cite{Rosen}, p.\ 325).
\begin{te}[Rubin]\label{T002}
If $\CT$ is a complete theory of linear orders, then $$I(\CT, \o )\in \{ 1, \c\}.$$
\end{te}
Following Fra\"{\i}ss\'{e} (see \cite{Fra}),  a structure $\Y$ of a relational signature $L$ is called {\it monomorphic}
iff,  for each positive integer $n$, all the $n$-element substructures of $\Y$ are isomorphic.
A complete theory $\CT \subset \Sent _L$ will be called a {\it monomorphic theory} iff each model $\Y$ of $\CT$
is monomorphic.
In Section \ref{S2} we show that the monomorphy of relational structures is an invariant of elementary equivalence and, moreover, that
a complete theory $\CT $ having infinite models is monomorphic iff it has a countable monomorphic model. In Section
\ref{S3} we prove the following statement, which is the main result of this paper.
\begin{te}\label{T030}
If $\CT$ is a complete monomorphic theory having infinite models, then $I(\CT ,\o)\in \{ 1 , \c \}$.
In addition,
$I(\CT ,\o)=  1$ iff some countable model of $\CT$ is simply definable by an $\o$-categorical linear order on its domain.
\end{te}
Section \ref{S4} contains some examples and comments. In particular, concerning the relationship between the class of monomorphic theories
and the aforementioned related classes of theories for which the Vaught conjecture was already confirmed,
we give a simple example of a complete monomorphic theory which is unstable, does not have definable Skolem functions
and has a countable model $\Y$ which is not bi-interpretable with any linear order. In addition, the linear orders which are related to $\Y$
(see Theorem \ref{T8067}) are neither o-minimal nor $\o$-categorical.

Concerning notation we note that throughout the paper
we assume that $L=\la R_i :i\in I\ra$ is a relational language, where  $\ar (R_i)=n_i\in \N$, for $i\in I$.
For $J\subset I$
by $L_J$ we will denote the reduction $\la R_i : i\in J\ra$ of $L$.
For convenience, $L_b$ will denote the binary language, $L_b=\la R \ra$, where $\ar (R)=2$,
and $L_\emptyset$ will be the empty language ($L_\emptyset$-formulas contain only the equality symbol).
For a theory $\CT \subset \Sent _L$, by $\Mod _{L}^{\CT }$ we denote the class of all models of $\CT$ and  by $\Mod _{L}^{\CT }(Y)$
the set of all models of $\CT$ with domain $Y$.
If $\Y \in \Mod _{L}(Y):=\Mod _{L}^{\emptyset }(Y)$ and $\emptyset \neq H\subset Y$, then $\BH$ will denote the corresponding substructure of $\Y$.

If $\X=\la X ,<\ra$ is a linear order, by $\X ^*$ we denote its reverse, $\la X ,< ^{-1}\ra $, and $\otp (\X )$ will denote the order type of $\X$.
For a set $X$ by $LO _X$ we denote the set of all linear orders on $X$.
\paragraph{Monomorphic structures}
For $n\in \N$, an $L$-structure $\Y=\la Y , \la R _i ^{\Y}:i\in I \ra \ra$ is called {\it $n$-monomorphic} iff all its substructures of size $n$ are isomorphic.
$\Y$ is said to be {\it monomorphic} iff it is $n$-monomorphic, for all $n\in \N$.

If $\Pa (\Y )$ denotes the set of all partial automorphisms of $\Y$, the structure $\Y$ is called
{\it chainable} if there is a linear order  $<$ on $Y$ such that $\Pa (\la Y,< \ra )\subset \Pa (\Y)$.

We will say that the structure $\Y$ is {\it simply definable in a linear order} 
iff there is a linear order $<$ on the set $Y$ such that
for each $i\in I$ there is a quantifier free $L_b$-formula $\f _i (v_0 ,\dots ,v_{n_i-1})$ which defines the relation $R _i ^{\Y}$
in the structure $\la Y ,< \ra$; that is,
$R _i ^{\Y}=D _{\la Y ,< \ra ,\f _i ,n_i}:=\{ \bar y\in Y^{n_i}: \la Y,<\ra  \models \f _i[\bar y]\}$.
Then we say that the linear order $<$ {\it chains} $\Y$.

Clearly, the oriented triangle is a finite monomorphic tournament which is not chainable.
But for infinite $L$-structures we have the following theorem,
proved by Fra\"{\i}ss\'{e} for finite languages (see \cite{Fra}) and for arbitrary languages by Pouzet \cite{Pou1}.
\begin{te}[Fra\"{\i}ss\'{e}]\label{T8067}
An infinite relational structure is monomorphic iff it
is chainable iff it is simply definable in a linear order.
\end{te}
So, by Theorem \ref{T8067}, an infinite $L$-structure $\Y$ is monomorphic iff the set
$\L _\Y:= \{ \la Y, \vartriangleleft \ra : \;\vartriangleleft \in LO_Y \mbox{ and }\la Y, \vartriangleleft \ra \mbox{ chains }\Y \}$
is non-empty and it is evident that $\X \in \L _\Y$ iff $\X ^* \in \L _\Y$. The following description of the set
$\L _\Y$ follows from Theorem 9 of  \cite{Gib},
which is a modification of similar results obtained independently by Frasnay in \cite{Fras}
and by Hodges, Lachlan and Shelah in \cite{Hodg1} (see also \cite{Fra}, p.\ 378, or \cite{Hodg}, p.\ 545).
\begin{te}[Gibson, Pouzet and Woodrow]\label{T8110}
If $\,\Y\in \Mod _L (Y)$ is an infinite monomorphic $L$-structure and $\X=\la Y ,<\ra\in \L _\Y$, then one of the following holds\\[-6mm]
\begin{itemize}\itemsep=-1mm
\item[{\sc (i)}] $\L _\Y=LO _Y$, that is, each linear order $\vartriangleleft$ on $Y$ chains $\Y$,
\item[{\sc (ii)}] $\L _\Y =\bigcup _{\X=\BI +\BF}\Big\{ \BF + \BI ,\, \BI ^* +\BF ^*\Big\}$,
\item[{\sc (iii)}] There are finite subsets $K$ and $H$ of $\,Y$ such that $\X =\BK +\BM +\BH\,$ and\\[2mm]
                   $ \L _\Y =\bigcup _{\vartriangleleft _K \in LO_K \atop \vartriangleleft _H \in LO_H}
                   \Big\{ \la K ,\vartriangleleft _K\ra + \BM +\la H ,\vartriangleleft _H \ra,
                   \la H ,\vartriangleleft _H\ra ^*+ \BM ^* +\la K ,\vartriangleleft _K \ra^*\Big\} $.
\end{itemize}
\end{te}
We note that the structures satisfying condition {\sc (i)} of the theorem are called {\it constant} by Fra\"{\i}ss\'{e}.
They are exactly the structures which are definable by the $L_\emptyset$-formulas (on their domain). The following simple fact will be used in Section \ref{S3}.
\begin{fac}\label{T8091}
If $\Y $ is a chainable $L$-structure 
and $\Z \cong \Y$, then $\otp [\L _\Y]=\otp [\L _{\Z}]$.
\end{fac}
\dok
Let $f\in \Iso (\Z , \Y)$, $\tau \in \otp [\L _\Y]$ and $\X =\la Y,< \ra\in \L _\Y$, where $\otp (\X )=\t$.
Then, clearly, $\X _1:=\la Z , f^{-1}[<]\ra \cong\X$ and $f\in \Iso (\X _1 ,\X )$; thus, $\otp (\X _1)=\tau$.
For $i\in I$ and $\bar z\in Z^{n_i}$ we have
$\bar z\in R_i^\Z$
iff $f\bar z\in R_i^\Y$
iff $\X \models \f _i [f\bar z]$
iff $\X _1\models \f _i [\bar z]$,
which gives $\X_1\in \L _{\Z}$. So, $\t =\otp (\X _1)\in \otp [\L _{\Z}]$ and we have proved that
$\otp [\L _\Y]\subset \otp [\L _{\Z}]$. The converse has a symmetric proof.
\kdok
\section{Monomorphic theories}\label{S2}
The notion of a monomorphic theory is established by the following theorem.
\begin{te}\label{T027}
If $\CT$ is a complete $L$-theory having infinite models, then the following conditions are equivalent:

(a) All models of $\;\CT$ are monomorphic,

(b)  $\CT$ has a monomorphic model,

(c)  $\CT$ has a countable monomorphic model.
\end{te}
The implication (a) $\Rightarrow$ (c) follows from  Claim \ref{T032}; the implication (c) $\Rightarrow$ (b) is trivial and
the implication (b) $\Rightarrow$ (a) follows from  Claim \ref{T028}(c) given below.
\begin{cla}\label{T029}
If $\BK$ is an $L$-structure of size $n\in \N$, then we have

(a) For each finite set $J\subset I$ there is an $L_J$-sentence $\p ^\BK _J$ such that for each $\Y \in \Mod _L$
\begin{equation}\label{EQ044}
\Y \models \p ^\BK _J \;\mbox{ iff }\;\, \forall H\in [Y]^n \;\;\BH |L_J \cong \BK |L_J ;
\end{equation}

(b) For the first-order theory $\CT ^\BK :=\{ \p ^\BK _J :J\in [I]^{<\o}\}$ and each $\Y \in \Mod _L$ we have
\begin{equation}\label{EQ045}
\Y \models \CT ^\BK \;\mbox{ iff }\;\; \forall H\in [Y]^n \;\;\BH \cong \BK .
\end{equation}
\end{cla}
\dok
(a) Let $K=\{ x_0 ,\dots ,x_{n-1}\}$ be an enumeration, let $\bar x=\la x_0 ,\dots ,x_{n-1}\ra$ and $J\in [I]^{<\o}$.
Writing $\bar v$ instead of $v_0 ,\dots ,v_{n-1}$, let $\Lit _{L_J}(\bar v)$ denote the set of all literals (atomic formulas and their negations) of $L_J$
with variables in the set $\{ v_0 ,\dots ,v_{n-1}\}$ and let
$\a ^\BK_J (\bar v) := \bigwedge \{ \eta \in \Lit _{L_J}(\bar v): \BK \models \eta [\bar x]\} $.
If $\Y$ is an $L$-structure, $\bar y \in Y^n$ an $n$-tuple and $H=\{ y_0 ,\dots ,y_{n-1}\}$,
then we have: $\Y \models \a ^\BK _J[\bar y]$ iff
$\{ \la x_k ,y_k\ra :k<n\}$ is an isomorphism from $\BK |L_J$ onto $\BH |L_J$.

It is easy to check that, in general, if $\pi\in \Sym (n)$ is a permutation and $\f _\pi (\bar v)$ is the formula obtained from a formula
$\f (\bar v)$ by replacement of $v_k$ by $v_{\pi (k)}$, for all $k<n$,
then $\Y \models \f _\pi [\bar y]$ iff $\Y \models \f [y_{\pi (0)}, \dots ,y_{\pi (n-1)} ]$. Thus
$\Y \models (\a ^\BK _J)_\pi [\bar y]$ iff
$p_\pi:=\{ \la x_k ,y_{\pi (k)}\ra :k<n\}$ is an isomorphism from $\BK |L_J$ onto $\BH |L_J$.
So, for the formula
$\f ^\BK_J (\bar v) := \bigvee _{\pi \in \Sym (n) }(\a ^\BK _J (\bar v))_\pi$ we have $\Y \models \f ^\BK _J[\bar y]$ iff $\BH |L_J \cong \BK |L_J$.
Now, for the sentence $\p ^\BK_J := \forall \bar v \; ( \bigwedge _{k<l<n}\neg v_k =v_l \Rightarrow \f ^\BK_J (\bar v))$
we have (\ref{EQ044}).

(b) Let $\Y \models \CT^\BK$ and suppose that $\BH \not \cong \BK$, for some $H=\{ y_0 ,\dots ,y_{n-1}\}\in [Y]^n$.
Then for each $\pi \in \Sym (n)$ we have $p_\pi \not\in \Iso (\BK ,\BH )$ and, since $p_\pi :K\rightarrow H$ is a bijection, there is $i_\pi\in I$
such that $p_\pi \not\in \Iso (\la K, R_{i_\pi}^{\BK}\ra ,\la H , R_{i_\pi}^{\BH }\ra )$.
Since $J:=\{ i_\pi : \pi \in \Sym (n)\}\in [I]^{<\omega}$ and $\Y \models \p ^\BK _J$, by (a)
there is $\pi _0\in \Sym (n)$ such that $p_{\pi _0}\in \Iso (\BK|L_J ,\BH |L_J )$, which implies that
$p_{\pi _0} \in \Iso (\la K, R_{i_{\pi _0}}^{\BK}\ra ,\la H , R_{i_{\pi _0}}^{\BH }\ra )$ and we have a contradiction.
The converse follows from (a).
\hfill $\Box$
\begin{cla}\label{T028}
If $\Y$ is a monomorphic $L$-structure and $K_n\in [Y]^n$, for $n\in \N$, then

(a) $\CT _M^\Y :=\bigcup _{n\in \N}\CT ^{\BK _n}\subset \Th (\Y )$;

(b) Each model $\Z$ of $\CT _M^\Y$ is monomorphic and $\Age (\Z)=\Age (\Y )$;

(c) If $\,\Z \models \Th (\Y )$, then $\Z $ is monomorphic and $\Age (\Z )=\Age (\Y )$.
\end{cla}
\dok
(a)  Since the structure $\Y$ is monomorphic, for $n\in \N$
we have: $\BH \cong \BK _n$, for all $H\in [Y]^{n}$; so, by Claim \ref{T029}(b), $\Y \models \CT ^{\BK _n}$.
Thus $\CT ^{\BK _n}\subset \Th (\Y )$, for all  $n\in \N$.

(b) If $\Z\models \CT _M^\Y$, then for each $n\in \N$ we have $\Z \models \CT ^{\BK _n}$ and, by  Claim  \ref{T029}(b),
$\BH \cong \BK _n$, for all $H\in [Z]^{n}$. Statement (c) follows from (a) and (b).
\hfill $\Box$
\begin{cla}\label{T032}
If $\CT$ is a complete monomorphic $L$-theory with infinite models and $|I|>\o$,
then $\CT$ has a countable model and there are a countable language $L_J \subset L$ and a complete monomorphic $L _J$-theory $\CT _J$
such that
\begin{equation}\label{EQ037}
\Big|\Mod _L ^\CT (\o )/\cong \Big| =\Big|\Mod ^{\CT _J}_{L_J}(\o)/\cong \Big|.
\end{equation}
\end{cla}
\dok
Let $\Y =\la Y, \la R_i^\Y :i\in I\ra\ra\in \Mod ^\CT _L$. By Theorem \ref{T8067}, there is a linear order $\la Y,<\ra$ which chains $\Y$.
Since there are countably many different relations (of all arities) defined by $L_b$-formulas in the structure $\la Y ,<\ra$,
there  is a partition $I=\bigcup _{j\in J}I_j$, where $|J|\leq \o$, such that, picking $i_j\in I_j$, for all $j\in J$,
we have $R_i^\Y = R_{i_j}^\Y$, for all $i\in I_j$.
So, for the $L$-sentences $\eta _{i,j}:=\forall \bar v \; (R_i (\bar v)  \Leftrightarrow R_{i_j} (\bar v))$, where $j\in J$ and $i\in I_j$, we have
$\CT _\eta := \bigcup _{j\in J}\{ \eta _{i,j} : i\in I_j \} \subset \Th _L (\Y )=\CT$.
Now, $L_J :=\la R_{i_j}:j\in J\ra \subset L$ and, using recursion,
to each $L$-formula $\f $ we adjoin an $L_J$-formula $\f _J$
in the following way:
$(v_k=v_l)_J :=v_k=v_l$;
$(R_i (v_{k_0}, \dots , v_{k_{n_i-1}}))_J := R_{i_j} (v_{k_0}, \dots , v_{k_{n_i -1}})$, for all $i\in I_j$;
$(\neg \f )_J :=\neg \f _J $; $(\f \land \p )_J :=\f _J \land \p _J$ and $(\forall v\; \f )_J :=\forall v\; \f _J $. A simple induction proves that
\begin{equation}\label{EQ036}
\forall \Z \in \Mod _L ^{\CT _\eta} \;\;\forall \f (\bar v)\in \Form _L \forall \bar z \in Z\;\;
\Big(\Z \models \f[\bar z] \Leftrightarrow \Z |L_J\models \f _J[\bar z] \Big) .
\end{equation}
In addition, for each $\Z _1 ,\Z _2 \in \Mod _L^{\CT _\eta}$ we have
\begin{equation}\label{EQ038}
\Z _1 \cong \Z _2 \Leftrightarrow \Z _1 |L_J\cong   \Z _2 |L_J \;\;\mbox{ and }\;\;\Z _1 \equiv _L \Z _2 \Leftrightarrow \Z _1 |L_J\equiv _{L_J}   \Z _2 |L_J .
\end{equation}
The first claim is true since $\Iso (\Z _1 ,\Z _2 )=\Iso (\Z _1 |L_J, \Z _2 |L_J )$.
For the second, suppose that $\Z _1 \equiv _L \Z _2$ and $\Z _1 |L_J \models \p$, where $\p \in \Sent _{L_J}$.
Then $\p \in \Sent _L$ and $\Z _1 \models \p$, which gives $\Z _2 \models \p$
so, by (\ref{EQ036}), $\Z _2 |L_J \models \p$, because $\p _J=\p$.
Conversely, suppose that $\Z _1 |L_J\equiv _{L_J} \Z _2 |L_J$ and $\Z _1\models \f$, where $\f \in \Sent _L$.
Then, by (\ref{EQ036}), $\Z _1 |L_J \models \f _J$ and, hence, $\Z _2 |L_J \models \f _J$, which by (\ref{EQ036}) gives $\Z _2\models \f$.

Let $\CT _J :=\Th _{L_J}(\Y |L_J)$. If $\Z \in \Mod _L^{\CT }$, that is, $\Z \equiv _L \Y$,
then by (\ref{EQ038}) we have   $\Z|L_J  \equiv _{L_J} \Y|L_J $, which means that $\Z|L_J \in \Mod _{L_J}^{\CT _J }$.
So we obtain the mapping $\Lambda : \Mod _L^{\CT } \rightarrow \Mod _{L_J}^{\CT _J }$, where $\Lambda (\Z )= \Z |L_J$, for all $\Z \in \Mod _L^{\CT }$,
which is an injection, because $\CT _\eta \subset\CT $.
If $\X \in \Mod _{L_J}^{\CT _J }$, then $\Z =\la X, \la R_i^{\Z } :i\in I\ra\ra\in \Mod _L^{\CT _\eta}$,
where $R_i^{\Z }=R_{i_j}^{\X }$, for $j\in J$ and $i\in I_j$. Now $\Z |L_J =\X \equiv _{L_J}\Y |L_J$ and, by (\ref{EQ038}), $\Z \equiv _{L}\Y $,
that is, $\Z\in \Mod _L^{\CT }$ and $\Lambda $ is a surjection.
Since the mapping $\Lambda$ preserves cardinalities of structures, we have
$\Lambda [\Mod ^\CT _L (\o)]=\Mod ^{\CT _J} _{L_J}(\o)$.
By the L\"{o}wenheim-Skolem theorem there is $\X \in \Mod ^{\CT _J}_{L_J}(\o)$ and
$\Lambda ^{-1}(\X)\in \Mod ^{\CT }_{L}(\o)$.

By (\ref{EQ038}), the mapping $\Lambda$ preserves the isomorphism relation and (\ref{EQ037}) is true.

Since the structure $\Y$ is monomorphic its reduct
$\Y |L_J$ is monomorphic as well and, by  Claim \ref{T028}(c), the theory $\CT _J$ is monomorphic.
\hfill $\Box$
\begin{rem}\label{R000}\rm
We list some comments and consequences of the claims given above.

1. {\it An $L$-structure is monomorphic iff all its finite reducts are monomorphic.}
For a proof of the non-trivial part, suppose that $|L|\geq \o$ and that  the reducts $\Y | L_J$, $J\in [I]^{<\o}$, are monomorphic.
Let us fix $n\in \N$ and $K\in [Y]^n$.
If $J\in [I]^{<\o}$, then $\BH | L_J \cong  \BK | L_J $, for all $H\in [Y]^n$; so, by  Claim \ref{T029}(a), $\Y \models \p ^\BK _J$.
Thus $\Y \models \CT ^\BK$ and, by  Claim \ref{T029}(b), $\BH \cong  \BK $, for all $H\in [Y]^n$;
so, 
$\Y$ is $n$-monomorphic.

2. {\it If $\CT$ is a complete monomorphic $L$-theory having infinite models, then all models of $\CT$ have the same age.}
This follows from Claim \ref{T028}(a). We note that, if $\Y$ is a linear order, then $\CT _M^\Y \varsubsetneq \Th (\Y )$,
because all linear orders have the same age.

3. Some of the sentences $\p ^{\BK }_J$ defined in Claim \ref{T029}(a) do not have infinite models;
for example, there is no $L_b$-structure of size $>3$ satisfying $\p ^{\BK }_{L_b}$, where $\BK$ is the 3-element circle tournament.

4. If $|L|<\o$, then, by Claim \ref{T029}, the first-order sentence $\p _n :=\bigvee _{\BK \in \Mod _L (n)} \p ^{\BK }_L$ says that an $L$-structure is $n$-monomorphic and the $L_{\o _1\o}$-sentence $\bigwedge _{n\in \N}\p _n $ says that an $L$-structure is  monomorphic.
But, by a theorem of Frasnay \cite{Fras} (see also \cite{Fra}, p.\ 359),
for each $n\in \N$ there is an integer $m\geq n$ such that for every infinite structure $\Y=\la Y, R^\Y\ra $ with one $n$-ary relation
we have:
$\Y$ is monomorphic
iff $\Y$ is ($\leq m$)-monomorphic
iff $\Y$ is $m$-monomorphic (Pouzet \cite{Pou}, see \cite{Fra}, p.\ 259).
So, the first-order sentence $\p _m :=\bigvee _{\r \subset m^n}\p ^{\la m,\r\ra}_{\la R\ra} $ says that an $n$-ary relation is monomorphic.
If $|L|=\k \geq\o$, then the monomorphy of $L$-structures  is expressed by the $L_{\k \o}$-sentence
$$\textstyle
\bigwedge _{n\in \N} \;\; \forall \bar v , \bar w \;\;
\Big[ \big(\bigwedge _{k<l<n}\neg v_k =v_l \land \bigwedge _{k<l<n}\neg w_k =w_l \big)\Rightarrow
$$
\vspace{-5mm}
$$\textstyle
\bigvee _{\pi \in \Sym (n) } \bigwedge _{i\in I}    \bigwedge _{\tau : n_i \rightarrow n}  \big( R_i (v_{\tau (0)}, \dots , v _{\tau (n_i)})
\Leftrightarrow  R_i (w_{\pi (\tau (0))}, \dots , w _{\pi (\tau (n_i))})\big)    \Big].
$$
\end{rem}
\section{Vaught's Conjecture}\label{S3}
In this section we prove Theorem \ref{T030}.
Let $\CT$ be a complete monomorphic $L$-theory having infinite models.
By  Claim \ref{T032}, w.l.o.g.\ we suppose that $|L|\leq \o$, which gives $\Mod_L^{\CT} (\o)\neq \emptyset$.
First we prove that
\begin{equation}\label{EQ026}
\Big|\Mod _{L}^{\CT }(\o )/ \!\cong \Big|\in \{ 1,\c\}.
\end{equation}
If $\Y _0=\la \o , \la R _i ^{\Y _0}:i\in I \ra \ra\in \Mod _{L}^{\CT }(\o )$,
then, by Theorem \ref{T027}, the structure $\Y _0$ is monomorphic and, by Theorem \ref{T8067},
there is a linear order $\X_0 \in \L _{\Y _0}\subset \Mod _{L_b}(\o )$
and for each $i\in I$ there is a quantifier free $L_b$-formula $\f _i (v_0 ,\dots ,v_{n_i-1})$ such that
\begin{equation}\label{EQ024}
\forall \bar x\in \o ^{n_i} \;\; \Big(\bar x\in R _i ^{\Y _0}\Leftrightarrow  \X _0 \models \f _i[\bar x]\Big).
\end{equation}
Clearly, $\CT =\Th _{L}(\Y _0)$. Let $\CT _{\X _0}$ denote the complete theory of $\X_0$, $\Th _{L_b}(\X _0)$.

Generally speaking, 
to each $L_b$-structure $\X\in \Mod _{L_b}(\o )$ we can adjoin the $L$-structure
$\Y _\X:=\la \o , \la R _i ^{\Y _\X}:i\in I \ra \ra \in \Mod _{L}(\o )$, where, for each $i\in I$,
the relation $R _i ^{\Y _\X}$ is defined in the structure $\X$ by the formula $\f_i$, that is,
\begin{equation}\label{EQ042}
\forall \bar x\in \o ^{n_i} \;\; \Big(\bar x\in R _i ^{\Y _\X }\Leftrightarrow  \X  \models \f _i[\bar x]\Big).
\end{equation}
\begin{cla}\label{T034}
For each structure $\Y _0\in \Mod _{L}^{\CT }(\o )$ and each linear order $\X_0 \in \L _{\Y _0}$, taking the formulas $\f _i$, $i\in I$, as above,
we have
\begin{itemize}
\item[\rm (a)] The mapping
$\Phi : \Mod _{L_b}(\o )\rightarrow \Mod _{L}(\o )$, defined by
$\Phi (\X )=\Y _\X$, for each $\X \in \Mod _{L_b}(\o )$,
preserves elementary equivalence and isomorphism.
Moreover, $\Iso (\X _1, \X _2)\subset \Iso (\Y _{\X _1}, \Y _{\X _2})$, for all  $\X _1, \X _2 \in \Mod _{L_b}(\o )$.
\item[\rm (b)] The mapping $\Psi : \Mod _{L_b}^{\CT _{\X _0}}(\o )/\!\cong  \;\rightarrow \Mod _{L}^{\CT }(\o )/\! \cong$, given by
$$
\Psi ([\X ])=[\Y _\X],
$$
for all $[\X ]\in \Mod _{L_b}^{\CT _{\X _0}}(\o )/\!\cong $,
is well defined.
\item[\rm (c)] For each linear order $\X\in \Mod _{L_b}^{\CT _{\X _0}}(\o )$ we have
$$
\Big|\Psi ^{-1}\Big[ \{[\Y _\X]\}\Big]\Big| \leq \Big|\otp \Big[\L _{\Y_{\X  }}\Big]\cap \otp \Big[\Mod _{L_b}^{\CT _{\X _0}}(\o )\Big]\Big|.
$$
\end{itemize}
\end{cla}
\dok
(a) By recursion on the construction of $L$-formulas
to each $L$-formula $\f (\bar v)$ we define an $L_b$-formula $\f ^*(\bar v)$ in the following way:
$(v_k=v_l)^* :=v_k=v_l$, $R_i (v_{k_0}, \dots , v_{k_{n_i-1}})^*:= \f _i (v_{k_0}, \dots , v_{k_{n_i-1}})$ (replacement of $v_j$ by $v_{k_j}$ in $\f_i$),
$(\neg \f)^* := \neg \f^*$, $(\f\land \p)^* := \f^*\land \p^*$ and $(\exists v_k\; \f)^* := \exists v_k\; \f^*$. A routine induction shows that, writing
$\bar v$ instead of  $v_0,\dots ,v_{n-1}$, we have
\begin{equation}\label{EQ025}
\forall \X \in \Mod _{L_b}(\o ) \;\;\forall \f (\bar v)\in \Form _L\;\forall \bar x \in \o ^n \;
\Big( \X \models \f ^*[\bar x]\Leftrightarrow \Y_\X \models \f[\bar x]\Big).
\end{equation}
Let $\X _1, \X _2 \in \Mod _{L_b}(\o )$. If $\X_1 \equiv \X_2$, then for an $L$-sentence $\f$ we have:
$\Y _{\X _1} \models \f$
iff  $\X _1\models \f ^*$ (by (\ref{EQ025}))
iff $\X _2 \models \f ^*$ (since $\X _1 \equiv \X_2$)
iff $\Y _{\X _2} \models \f$ (by (\ref{EQ025}) again). So, $\Y _{\X _1} \equiv \Y _{\X_2} $ and the mapping $\Phi $
preserves elementary equivalence.

If $f:\X_1 \rightarrow  \X _2$ is an isomorphism,
then by (\ref{EQ042}) and since  isomorphisms preserve all formulas in both directions, for each  $i\in I$ and $\bar x\in \o ^{n_i}$ we have:
$\bar x\in R_i ^{\Y_{\X _1}}$
iff $\X _1 \models \f _i [\bar x]$
iff $\X _2 \models \f _i [f\bar x]$
iff $f \bar x\in R_i ^{\Y_{\X _2}}$.
Thus $f\in \Iso (\Y _{\X_1}, \Y _{ \X _2})$.

(b) For $\X \in \Mod _{L_b}^{\CT _{\X _0}}(\o )$ we have $\X \equiv \X_0$, which, by (a),  (\ref{EQ024}) and (\ref{EQ042}),
implies that $\Phi (\X )=\Y _\X \equiv \Y _{\X_0}=\Y _0$.
So, since $\Y_0\models \CT$, we have $\Phi (\X )\in \Mod _{L}^{\CT }(\o )$ and, thus,
\begin{equation}\label{EQ043}
\Phi [ \Mod _{L_b}^{\CT _{\X _0}}(\o )] \subset \Mod _{L}^{\CT }(\o ).
\end{equation}
Assuming that $\X_1, \X _2 \in \Mod _{L_b}^{\CT _{\X _0}}(\o )$ and $\X _1\cong \X _2$, by (a) we have $\Y_{\X _1}\cong \Y_{\X _2}$,
that is $[\Y_{\X _1}]=[\Y_{\X _2}]$. So, the mapping $\Psi $ is well defined.

(c) By (\ref{EQ043}) and (b), $\Y _\X \in \Mod _{L}^{\CT }(\o )$ and
$\Psi ^{-1}[\{ [\Y _\X]\}]\subset \Mod _{L_b}^{\CT _{\X _0}}(\o )/\cong $.
We show that the correspondence $[\X _1]\mapsto \otp (\X _1)$ is an injection from the set  $\Psi ^{-1}[\{[\Y _\X]\}]$
to the set of order types $\;\otp [\L _{\Y_{\X  }}]\cap \otp [\Mod _{L_b}^{\CT _{\X _0}}(\o )]$.

For $[\X _1]\in \Psi ^{-1}[\{[\Y _\X]\}]$ we have $[\Y _{\X _1}]=\Psi ([\X _1])=[\Y _\X] $, that is, $\Y _{\X _1}\cong\Y _\X $ and,
since $\X _1\in \L _{\Y_{\X _1 }}$, by Fact \ref{T8091} we have
$\otp (\X _1)\in \otp[\L _{\Y_{\X _1 }}]=\otp[\L _{\Y_\X }]$.
Also we have $\X_1\in \Mod _{L_b}^{\CT _{\X _0}}(\o )$ and, hence,
$\otp(\X_1)\in \otp[\Mod _{L_b}^{\CT _{\X _0}}(\o )]$.

In addition, if $[\X _1],[\X _2]\in \Psi ^{-1}[\{[\Y _\X]\}]$
and $[\X _1]\neq[\X _2]$, then $\X _1\not\cong\X _2$, and, hence, $\otp (\X _1)\neq \otp (\X _2)$.
Thus $[\X _1]\mapsto \otp (\X _1)$ is an injection indeed.
\kdok
\begin{cla}\label{T019'}
If some structure $\Y \in \Mod _{L}^{\CT }(\o )$ is chained by an $\o$-categorical linear order,
then $\Y $ is an $\o$-categorical $L$-structure.
\end{cla}
\dok
Let $\Y \in \Mod _{L}^{\CT }(\o )$ and let $\X\in \L _\Y$ be an $\o$-categorical linear order.

By the theorem of Engeler, Ryll-Nardzewski and Svenonius (see \cite{Hodg}, p.\ 341),
the automorphism group of $\X $ is oligomorphic; that is, for each $n\in \N$ we have
$|\o ^n /\!\sim _{\X ,n}|<\o$, where $\bar x \sim_{\X ,n} \bar y$ iff $f\bar x =\bar y$, for some $f\in \Aut (\X )$.

As in Claim \ref{T034}(a) we prove that $\Aut (\X )\subset \Aut (\Y )$, which implies that for $n\in \N$
and each $\bar x , \bar y \in \o ^n$ we have $\bar x \sim_{\X ,n} \bar y \Rightarrow \bar x \sim_{\Y ,n} \bar y$.
Thus $|\o ^n /\!\sim _{\Y ,n}| \leq |\o ^n /\!\sim _{\X ,n}|<\o$, for all $n\in \N$, and, since $|L|\leq \o$, using the same theorem we conclude that
$\Y $ is an $\o$-categorical $L$-structure.
\kdok
Now we prove (\ref{EQ026}) distinguishing the following cases and subcases.\\[-2mm]

\noindent
{\bf Case A:} {\it Some structure $\Y \!\in \Mod _{L}^{\CT }(\o )$ is chained by an $\o$-categorical linear order.}
Then, by Claim \ref{T019'}, the structure $\Y $ is $\o$-categorical, that is,
$|\Mod _{L}^\CT (\o )/ \!\cong | = 1$.\\[-2mm]

\noindent
{\bf Case B:} {\it The set $\bigcup _{\Y \!\in \Mod _{L}^{\CT }(\o )}\L _\Y$ does not contain $\o$-categorical linear orders.}
Then, by Theorem \ref{T002} and since there is no structure $\Y \in \Mod _{L}^{\CT }(\o )$ which is chained by each linear order on $\o$, we have
\begin{equation} \label{EQ027}\textstyle
\forall \Y \!\in \Mod _{L}^{\CT }(\o ) \;\;\forall \X \in \L _\Y\;\;
\Big|\Mod _{L_b}^{\CT _{\X }}(\o )/\!\cong\Big|=\c ,
\end{equation}
\begin{equation} \label{EQ040}
\forall \Y \in \Mod _{L}^{\CT }(\o ) \;\; \L _\Y \neq LO _\o .
\end{equation}
We prove that $|\Mod ^\CT _L (\o )/\!\cong|=\c$, distinguishing the following two subcases.\\[-2mm]

\noindent
{\bf Subcase B1:} {\it For some $\Y _0\in \Mod _{L}^{\CT }(\o )$ there is a linear order $\X_0 \in \L _{\Y _0}$ having at least one end-point.}

Then by (\ref{EQ027}) we have $|\Mod _{L_b}^{\CT _{\X _0}}(\o )/\!\cong|=\c$
and, by Claim \ref{T034}(b), for a proof that $|\Mod ^\CT _L (\o )/\!\cong|=\c$
it is sufficient to show that the mapping $\Psi$ is at-most-countable-to-one (see Example \ref{EX011}). This will follow from the following claim and
Claim \ref{T034}(c).
\begin{cla}\label{T023}
$\Big|\otp [\L _{\Y_{\X  }}]\cap \otp [\Mod _{L_b}^{\CT _{\X _0}}(\o )]\Big|\leq \o$, for all $\X\in \Mod _{L_b}^{\CT _{\X _0}}(\o )$.
\end{cla}
\dok
Let $\X\in \Mod _{L_b}^{\CT _{\X _0}}(\o )$ and $\t :=\otp (\X)$.
First, if the set $\L _{\Y_{\X }}$ satisfies {\sc (iii)} of Theorem \ref{T8110}, then we have $\otp[\L _{\Y_{\X }}]=\{ \t, \t^*\}$
and the claim is proved.

Otherwise, by (\ref{EQ040}) and Theorem \ref{T8110} we have
$\L _{\Y_\X} =\bigcup _{\X=\BI +\BF}\{ \BF + \BI ,\, \BI ^* +\BF ^*\}$.
Let $\X=\BI +\BF$, where $\BI,\BF\neq\emptyset$.

If $\BI$ has a largest element, say $x$; then $\BI =(-\infty ,x]_\X$ and $\BF =(x,\infty )_\X$;
   so, $\t _x :=\otp (\BF + \BI )=\otp((x,\infty )_\X +(-\infty ,x]_\X)$ and $\t _x^* :=\otp (\BI ^* +\BF ^*)=\otp ((-\infty ,x]_\X^* +(x,\infty )_\X^*)$
   are the corresponding elements of $\otp [\L _{\Y_{\X  }}]$.

If $\BI$ does not have a largest element and $\BF$ has a smallest element, say $x$;
   then $\BI =(-\infty ,x)_\X $ and $\BF =[x,\infty )_\X$;
   so, $\s _x :=\otp (\BF + \BI )=\otp([x,\infty )_\X +(-\infty ,x)_\X )$ and $\s _x^* :=\otp (\BI ^* +\BF ^*)=\otp ((-\infty ,x)_\X^* +[x,\infty )_\X^*)$
   are the corresponding elements of $\otp [\L _{\Y_{\X  }}]$.

If $\BI$ does not have a largest element and $\BF$ does not have a smallest element, then $\{ \BI,\BF\}$
   is a gap in $\X$ so $\BF +\BI$ and $\BI ^* +\BF ^*$ are linear orders without end points.
   By our assumption, the linear order $\X_0$ has at least one end point and, since this is a first-order property,
   we have $\BF +\BI ,\BI ^* +\BF ^*\not\equiv \X_0 $.
   Consequently we have $\otp (\BF + \BI ), \otp (\BI ^* +\BF ^*)\not\in \otp [\Mod _{L_b}^{\CT _{\X _0}}(\o )]$.

Thus, $\otp [\L _{\Y_{\X  }}]\cap \otp [\Mod _{L_b}^{\CT _{\X _0}}(\o )]\subset
\{ \t, \t^*\} \cup \bigcup_{x\in \o}\{\t _x ,\t _x ^* ,\s _x , \s _x ^* \}=:\Theta$
and, since $\Theta$ is a countable collection of order types, the claim is proved.
\kdok

\noindent
{\bf Subcase B2:} {\it Each $\X \in \bigcup _{\Y \in \Mod _{L}^{\CT }(\o )} \L _{\Y }$ is a linear order without end points.}

Then, by (\ref{EQ040}) and Theorem \ref{T8110}, for each $\Y \in \Mod _{L}^{\CT }(\o )$  the set
$\L _{\Y }$ is of the form {\sc (iii)} or {\sc (ii)}.
Let us fix arbitrary $\Y _0 \in \Mod _{L}^{\CT }(\o )$ and $\X _0 \in \L _{\Y _0 }$. By (\ref{EQ027})
we have $|\Mod _{L_b}^{\CT _{\X _0} }(\o )/\cong|=\c$ and  Claim \ref{T034} is true for $\Y _0 $ and $\X _0 $.

Suppose that there is $\X \in \Mod _{L_b}^{\CT _{\X _0} }(\o )$ such that $\L _{\Y_{\X } }$ is given by {\sc (ii)},
that is $\L _{\Y_{\X } } =\bigcup _{\X =\BI +\BF}\{ \BF + \BI ,\, \BI ^* +\BF ^*\}$.
Then, since $\X $ is a linear order without end points,
picking an arbitrary $x\in \o$
we have $\X =(-\infty ,x)_{\X }+ [x, \infty)_{\X }$; thus,
the linear order $[x, \infty)_{\X }+ (-\infty ,x)_{\X }$ chains $\Y_{\X }\in \Mod _{L}^{\CT }(\o )$ and has a minimum, which contradicts the assumption
of Subcase B2.

So, for each $\X \in \Mod _{L_b}^{\CT _{\X  _0}}(\o )$ the set $\L _{\Y_{\X } }$ is given by {\sc (iii)}, for some $n\in \o$.
In addition, since each element of $\L _{\Y _\X }$ is a linear order without end points, we have $K=H=\emptyset$ and, hence, $\L _{\Y _\X }=\{ \X ,\X ^*\}$,
which gives $|\otp [\L _{\Y_{\X } }]|\leq 2$. Now, as above, we obtain $|\Mod _{L}^\CT (\o )/ \!\cong | = \c$ and the proof of (\ref{EQ026}) is finished.

If some countable model of $\CT$ is chained by an $\o$-categorical linear order, then by Claim \ref{T019'} we have $|\Mod _{L}^\CT (\o )/ \!\cong | =1$.
Otherwise we have Case B and, as above,
$|\Mod _{L}^\CT (\o )/ \!\cong | = \c$. Thus $|\Mod _{L}^\CT (\o )/ \!\cong | =1$ iff some $\Y \in \Mod _{L}^\CT (\o )$
is chained by some $\o$-categorical linear order.
\kdok
\section{Examples}\label{S4}
\begin{ex}\label{EX010}\rm
{\it A complete monomorphic theory which is unstable, does not have definable Skolem functions
and has a countable model which is not bi-interpretable with any linear order.}

Let $Q$ denote the set of rational numbers and $Z$ the set of integers.
Let $\X =\la X,<\ra$ be the linear order $\sum _{q\in \Q}L_q$, where
$L_q =\{ q\}$, for $q\in Q\setminus Z$; $L_q =\{ n,n'\}$, for $q=n\in Z$, and $n<n'\not\in Q$.
The $L_b$-formula
$$
\f _{\mathrm{betw}}(v_0 ,v_1 , v_2):=v_0 <v_1 < v_2 \lor v_2 <v_1 < v_0
$$
defines the betweness relation $D_{\f _{\mathrm{betw}}}=\{ \bar x\in X^3 : \X \models \f _{\mathrm{betw}}[\bar x]\}$
on the set $X$ and $\Y =\la X , D_{\f _{\mathrm{betw}}}\ra$ is a monomorphic $L_3$-structure, where $L_3 =\la S\ra$ and $\ar (S)=3$.

In order to prove that the theory $\Th (\Y )$ is unstable we show that it has the order property (see \cite{Hodg}, p.\ 307).
Let $\f (v_0,v_1,v_2, w_0,w_1,w_2 )=\f (\bar v ,\bar w)$ be the $L_3$-formula
$$
S(v_1 ,v_2 ,v_0) \land S(v_1 ,v_2 ,w_0) \land S(v_2 ,v_0 ,w_0)
$$
and let $\bar x _k =\la k,-2,-1\ra \in X^3$, for $k\in \o$. Then for $k,l\in \o$ we have
$\Y \models \f [\bar x _k ,\bar x _l]$,
iff $\Y \models \f [ k,-2,-1 , l,-2,-1]$,
iff $\Y \models S[-2 ,-1 ,k]$ and $\Y \models S[-2 ,-1 ,l]$ and $\Y \models S[-1 ,k ,l]$,
iff $k,l> -1$ and $k$ is between $-1$ and $l$,
iff $k<l$.

Suppose that $\Th (\Y )$ has definable Skolem functions\footnote{that is, (see \cite{Hodg}, p.\ 91)

$\forall \f (\bar v,w)\in \Form _{L_3} \;\;
\exists \p _\f (\bar v,w)\in \Form _{L_3} \;\;
\forall \X \in \Mod _L ^{\Th (\Y )}\;\;
\forall \bar x\in X^n \;\;
\forall y\in X $

$\big[
\X \models \f [\bar x ,y] \Rightarrow
\exists z\in X \;\;
\big(
\X \models \p _\f [\bar x , z]\land \f [\bar x , z] \land
\forall t\in X\setminus \{ z\}\;\;
\X \models \neg \p _\f [\bar x , t]
\big)
\big].
$
}
and that $\p _S (v_0,v_1,v_2 )$ is an $L_3$-formula corresponding to the atomic $L_3$-formula $S (v_0,v_1,v_2 )$.
Then for $\Y$ and $0' , 1 \in X$, since $\Y \models S [0',  \frac{1}{2} ,1]$,
there would be a $z\in X $ such that:
(i) $\Y \models S [0',  z ,1]$, which means that $0' <z <1$,
(ii) $\Y \models \p _S [0', z, 1]$,
(iii) $\Y \models \neg \p _S [0', t , 1]$, for all $t\in X\setminus \{ z\}$.
Let $t\in (0' ,1)_\X \setminus \{ z\}$ and let $f\in \Aut (\X )$, where $f(0')=0'$, $f(1)=1$ and $f(z)=t$.
Then $f\in \Aut (\Y )$ as well and by (ii) we have $\Y \models \p _S [0', t, 1]$, which contradicts (iii).

Suppose that there is a linear order $\BL$ such that the  structures $\Y$ and $\BL$ are bi-interpretable.
Then we would have $\Aut (\Y )\cong \Aut (\BL)$ (see \cite{Hodg}, p.\ 226).
Let $f:X\rightarrow X$ be the strictly $<$-decreasing bijection defined by $f(q)=-q$, for $q\in Q\setminus Z$,
and $f(n)=-n'$ and $f(n')=-n$, for $n\in Z$; (thus $f(0)=0'$ and $f(0')=0$).
It is easy to see that $f$ reverses the linear order $<$ and, hence, $f\in \Aut (\Y )$.
It is evident that  $f\circ f =\id _X$ and $f\neq \id _X$; thus, since
$\Aut (\BL)$ does not contain non-zero elements of order 2, we have a contradiction.
(Suppose that $g\in \Aut (\BL)$, $g\circ g =\id _L $ and $g\neq \id _L $. Then $g(x)\neq x$,
for some $x\in L$, say $g(x)<x$; but then $x=g(g(x))<g(x)$ and we have a contradiction.)
Moreover, no linear order is a retraction of $\Y$
(otherwise we would have an injective homomorphism $h:\Aut (\Y )\hookrightarrow \Aut (\BL)$ and, hence, $h(f)\circ h(f)=\id _L$ and $h(f)\neq\id _L$).

Suppose that some linear order $\vartriangleleft$ is definable in the structure $\Y$.
Then there is an $L_3$-formula $\p (v_0,v_1)$ such that
\begin{equation}\label{EQ041}
\forall x_0 ,x_1 \in X \;\;\Big(  x_0 \vartriangleleft x_1 \Leftrightarrow \Y \models \p [x_0 ,x_1]\Big).
\end{equation}
Let $f\in \Aut (\Y )$ be the function from the previous paragraph.
Now, if $0\vartriangleleft 0'$, then by (\ref{EQ041}) we have $\Y \models \p [0 ,0']$ and, hence, $\Y \models \p [f(0) ,f(0')]$,
that is $\Y \models \p [0' ,0]$ so, by (\ref{EQ041}), $0'\vartriangleleft 0$, which  is impossible.
In a similar way we show that $0'\vartriangleleft 0$ is impossible and we have a contradiction.
In particular, the structures $\X$ and $\Y$ are not bi-definable.

We recall that o-minimal structures (widely investigated by Pillay and Steinhorn, see \cite{Pil}) are
linearly ordered structures in which every parametrically definable subset of the domain
is a finite Boolean combination of intervals. Since the set $Z$ is definable in the linear order $\X$ as the set
of elements having an immediate successor,  the linear order $\X$ is not o-minimal.

In addition, $\X$ is not an $\o$-categorical linear order, because for each $n\in \o$
there are $x,y\in X$ such that the interval $(x,y)_\X$ contains exactly $n$ elements having an immediate successor; thus $|S_2 (\Th (\X))|\geq \o$.
(See also Rosenstein's characterization of $\o$-categorical theories of linear orders, \cite{Rosen}, p.\ 299.)
Similarly, $\Th (\Y)$ is not an $\o$-categorical theory and by Theorem \ref{T030} we have $I(\Th (\Y ) ,\o)= \c $.
\end{ex}
\begin{ex}\label{EX011}\rm
{\it The mapping $\Psi$ from Claim \ref{T034} must not be finite-to-one.}

Let $\X_0$ be a linear order of the type $\o +\o ^*$. Then the countable models of $\Th (\X _0)$
are of the form $\X _\tau :=\o + \zeta \tau +\o ^*$, where $\tau$ is a countable (or the empty) order-type and $\zeta$ the order type of the integers.
The $L_b$-formula
$$
\f _{\mathrm{cyclic}}(v_0 ,v_1 , v_2):=v_0 <v_1 < v_2 \lor v_1 <v_2 < v_0 \lor v_2 <v_0 < v_1
$$
defines the  cyclic-order relation on the linear order $\X _\o \in \Mod _{L_b}^{\Th (\X _0)}(\o )$, say $\X _\o=\la \o ,<_{\X _\o} \ra\cong\o + \zeta \o +\o ^*$.
So, defining $D_{\f _{\mathrm{cyclic}}}=\{ \bar x\in \o^3 : \X _\o \models \f _{\mathrm{cyclic}}[\bar x]\}$
we obtain a monomorphic $L_3$-structure $\Phi (\X _\o ) =\Y _{\X _\o} =\la \o , D_{\f _{\mathrm{cyclic}}}\ra \in \Mod _{L_3}^{\Th (\Y _{\X _0})}$
(where $L_3 =\la S\ra$ and $\ar (S)=3$ again).
It is well known (see, for example, \cite{Fra}) that $\L _{\Y _{\X _\o}}=\bigcup _{\X _\o =\BI +\BF }\{ \BF +\BI , \BI ^* +\BF ^* \}$ and using that equality
we easily check that
$$
\otp\Big[\L _{\Y _{\X _\o}}\Big]=\bigcup _{n\in \o}\Big\{ \o + \zeta (\o +n) +\o ^* , \o + \zeta (n+\o ^*) +\o ^*, \zeta (\o +n) , \zeta (n+\o ^*)\Big\} .
$$
So, since the models of the theory $\Th (\X _0)$ have end points, we obtain
$$
\otp\Big[\L _{\Y _{\X _\o}}\Big] \cap \;\otp \Big[\Mod _{L_b}^{\Th (\X _0)}(\o )\Big]
=\bigcup _{n\in \o}\Big\{ \o + \zeta (\o +n) +\o ^* , \o + \zeta (n+\o ^*) +\o ^*\Big\}
$$
and, hence, $\Big|\Psi ^{-1}\Big[[\Y _{\X _\o} ]\Big]\Big|=\o$.
\end{ex}

\noindent
{\bf Acknowledgments.}
This research was supported by the Ministry of Education and Science of the Republic of Serbia (Project 174006).

\end{document}